\documentclass[12pt,preprint,spbasic]{aastex}
\usepackage{spr-astr-addons}
\usepackage{aas_macros}
\usepackage{natbib}
\usepackage{url}
\newcommand{\amc}{{(1-\mu) }}
\newcommand{\zabs}{{\Bigl(1+\frac{5A_2}{2{r^2_2}_*}\Bigr)}}
\shorttitle{Linear Stability of  Equilibrium Points ---- Chermnykh's Problem}

\shortauthors{Badam Singh Kushvah}

\begin{document}
\title{Linear Stability of  Equilibrium Points in the Generalized Photogravitational Chermnykh's Problem}

\author{Badam Singh Kushvah\altaffilmark{1}}
\affil{Gwalior Engineering College  Airport  Road, Maharajpura,
Gwalior (M.P.)-474015,INDIA}\email{bskush@gmail.com}

 \altaffiltext{1}{Present
address for correspondence: 139/I, Anupam Nagar Ext-2, Opp. Jiwaji University, City Center,
Gwalior-474011} 

%

\begin{abstract}
The  equilibrium points and their linear stability has been  discussed in the generalized photogravitational Chermnykh's problem. The bigger primary is being  considered as   a source of radiation and small primary as an oblate spheroid. The effect of radiation pressure   has been discussed numerically.  The collinear points are linearly unstable and triangular points  are  stable in the sense of Lyapunov stability provided  $\mu< \mu_{Routh}=0.0385201$. The effect of  gravitational potential from the belt is also examined.  The mathematical properties of this system are different from the classical restricted three body problem. 
\end{abstract}



%
%
%
\section{Introduction}
\label{intro}
The Chermnykh's problem is new kind of restricted three body problem which was first time studied by \citet*{Chermnykh1987VeLen}. Recently many authors studied this problem such as \citet*{JiangYeh2004AJ,JiangYeh2004MNRAS,JiangYeh2004RMxAC}  considered the influence from the belt for planetary systems and found that the probability  to have equilibrium points around the inner part of the belt is larger than the one near the outerpart.  \citet*{Papadakis2005Ap&SS} examined the motion around the triangular equilibrium points of the restricted three-body problem under angular velocity variation.
\citet*{YehJiang2006Ap&SS} studied a Chermnykh-Like problem in which the mass parameter $\mu$ is set to be $0.5$. \citet*{JiangYeh2006Ap&SS} found the equilibrium points in the Chermnykh-Line problem when an addition gravitational potential from the belt is included.
The solar radiation pressure force $F_p$ is exactly apposite to the gravitational attraction force $F_g$ and change with the distance by the same law it is possible to consider that the result of action of this force will lead to reducing the effective mass of the  Sun or particle. It is acceptable to speak about a reduced mass of the particle as the effect of reducing its mass depends on the properties of the particle itself.
 
\citet*{KushvahBR2006} examined the linear stability of triangular
equilibrium points in the generalized photogravitational restricted
three body problem with Poynting-Robertson drag, $L_4$ and $L_5$
points became unstable due to P-R drag which is very remarkable and
important, where as they are linearly stable in classical problem when $0<\mu<\mu_{Routh}=0.0385201$.  \citet*{KushvahetalHr2007,KushvahetalNON2007,KushvahetalNr2007} examined normalization of Hamiltonian they have also studied the nonlinear stability of triangular equilibrium points in the generalized photogravitational restricted three body problem with Poynting-Robertson drag, they have found  that the triangular points are stable in the nonlinear sense except three critical mass ratios at which KAM theorem fails. \citet*{PapadakisKanavos2007Ap&SS}  given numerical exploration of Chermnykh's problem, in which the equilibrium points and zero velocity curves studied numerically also the non-linear stability for the triangular Lagrangian points are computed numerically for the Earth-Moon and Sun-Jupiter mass distribution when the angular velocity varies. The stability of triangle libration points in generalized restricted circular three-body problem has been studied by \citet*{Beletsky2008CosRe}. \citet*{MKDas2008Ap&SS} examined the stability of location of various equilibrium points of a passive micron size particle in the field of radiating binary stellar system within the framework of circular restricted three body problem.
 
 In this paper we have obtained the equations of motion, the position of equilibrium points and their linear stability in the generalized photogravitaional Chermnykh's problem.  The effect of radiation pressure, oblateness, and  gravitational potential from the belt  has been  examined  analytically and numerically. We have seen that the collinear equilibrium  points are linearly unstable while the triangular points are conditionally stable.
\section{Equations of Motion and Zero Velocity Curves}
\label{sec:eqmot} 
We consider the barycentric rotating co-ordinate system $Oxyz$ relative to inertial system with angular velocity $\omega$ and common $z$--axis.  We have taken line joining the primaries as $x$--axis. Let $m_1, m_2$ be the masses of bigger primary(Sun)  and smaller primary(Earth) respectively. Let  $Ox$, $Oy$  in the equatorial plane of smaller primary  and $Oz$ coinciding with the polar axis of $m_{2}$. Let $r_{e}$, $r_{p}$ be the equatorial and polar radii of $m_{2}$ respectively,  $r$ be the distance between primaries.  Let infinitesimal mass $m$ be placed at the  point $P(x,y,0)$. We take units such that sum of the masses and distance between primaries is  unity, the unit of time i.e. time period of $m_{1}$ about $m_{2}$  consists of $2\pi$ units such that the Gaussian constant of gravitational $\Bbbk^{2}=1$. Then perturbed mean motion $n$ of the primaries is given by $n^{2}=1+\frac{3A_{2}}{2}$, where $A_{2}=\frac{r^{2}_{e}-r^{2}_{p}}{5r^{2}}$ is oblateness coefficient of $m_{2}$.
Let $\mu=\frac{m_{2}}{m_{1}+m_{2}}$ then $1-\mu=\frac{m_{1}}{m_{1}+m_{2}}$ with $m_{1}>m_{2}$, where $\mu$ is mass parameter. Then coordinates of $m_{1}$ and $m_{2}$ are  $(-\mu,0)$ and $(1-\mu,0)$ respectively. 
In the above mentioned reference system  we determine  the equations of motion of the infinitesimal mass particle in $x y$-plane  as \citet*{kushvah-2008}.
\begin{eqnarray}
\ddot{x}-2n\dot{y}&=&U_{x} ,\label{eq:ux}\\
\ddot{y}+2n\dot{x}&=&U_{y} \label{eq:uy} \end{eqnarray}
where
\begin{eqnarray*}
&&U_x= n^{2}x-\frac{(1-\mu)q_1(x+\mu)}{r^3_1}-\frac{\mu(x+\mu-1)}{r^3_2}\\&&-\frac{3}{2}\frac{\mu{A_2}(x+\mu-1)}{r^5_2} \nonumber\\
&&U_y=n^{2}y
-\frac{(1-\mu)q_{1}{y}}{r^3_1}
-\frac{\mu{y}}{r^3_2}-\frac{3}{2}\frac{\mu{A_2}y}{r^5_2} \nonumber\end{eqnarray*}
where
\begin{eqnarray}
&&U=\frac{n^2(x^2+y^2)}{2}+\frac{(1-\mu)q_1}{r_1}+\frac{\mu}{r_2}+\frac{\mu
 A_2}{2r_2^3}\label{eq:FF}
 \end{eqnarray}
 $q_1=1-\frac{F_p}{F_g}$ is a mass reduction factor expressed in terms of the particle radius $\mathbf{a}$, density $\rho$ radiation pressure efficiency factor $\chi$ (in C.G.S. system): \( q_1=1-\frac{5.6\times{10^{-5}}}{\mathbf{a}\rho}\chi
 \leq1 \) . The assumption $q_1=constant$ is equivalent to neglecting fluctuations in the beam of solar radiation and the effect of the planets shadow.
\subsection{\citet*{MiyamotoNagai1975PASJ}  Profile Model}
\label{subsect:modelA} In this model we introduce the  potential of belt as:
\begin{equation}
 V(r,z)=-\frac{M_b}{\sqrt{r^2+\left(\mathbf{a}+\sqrt{z^2+\mathbf{b}^2}\right)^2}}\label{eq:VRZ}
\end{equation}
where $M_b$ is the total mass of the belt and $r^2=x^2+y^2$, $\mathbf{a,b}$ are parameters which determine the density profile of the belt. The parameter $\mathbf{a}$ controls the flatness of the profile and can be called  \lq\lq flatness parameter\rq\rq. The parameter $\mathbf{b}$ controls the size of the core of the density profile and can be called \lq\lq core parameter\rq\rq. When $\mathbf{a}=\mathbf{b}=0$ the potential equals to the one by a points mass. In genrel the density distribution  corresponding to the above $V(r,z)$ in (~\ref{eq:VRZ}) is as in \citet*{MiyamotoNagai1975PASJ}
\begin{eqnarray}
 \rho (r,z)=\frac{\mathbf{b}^2M_b\left[\mathbf{a}r^2 +\left(\mathbf{a}+3N\right)\right]\left( \mathbf{a}+N\right)^2}{N^3\left[r^2+\left(\mathbf{a}+N\right)^2\right]^{5/2}}
\end{eqnarray}
 where $N=\sqrt{z^2+\mathbf{b}^2}$,\ $T=\mathbf{a}+\mathbf{b}$, $z=0$.
Then we obtained \begin{equation}
 V(r,0)=-\frac{M_b}{\sqrt{r^2+T^2}}\label{eq:Vr0}
\end{equation}
and  \(V_x=\frac{M_b x}{\left(r^2+T^2\right)^{3/2}},\ V_y=\frac{M_b y}{\left(r^2+T^2\right)^{3/2}} \).
 Now we consider  only the orbits on the $x-y$ plane, then the equations of motion are modified  by using (~\ref{eq:ux},\ref{eq:uy}) in  the following form:
\begin{eqnarray}
\ddot{x}-2n\dot{y}&=&U_{x}-V_x=\Omega_x ,\label{eq:Omegax}\\
\ddot{y}+2n\dot{x}&=&U_{y}-V_y=\Omega_y\label{eq:Omegay}
 \end{eqnarray}
where
\begin{eqnarray*}
&&\Omega_x= n^{2}x-\frac{(1-\mu)q_1(x+\mu)}{r^3_1}-\frac{\mu(x+\mu-1)}{r^3_2}\\&&-\frac{3}{2}\frac{\mu{A_2}(x+\mu-1)}{r^5_2}-\frac{M_b x}{\left(r^2+T^2\right)^{3/2}} \nonumber\\
&&\Omega_y=n^{2}y
-\frac{(1-\mu)q_{1}{y}}{r^3_1}
-\frac{\mu{y}}{r^3_2}\nonumber\\&&-\frac{3}{2}\frac{\mu{A_2}y}{r^5_2}-\frac{M_b y}{\left(r^2+T^2\right)^{3/2}}\end{eqnarray*}

\begin{eqnarray}
&&\Omega=\frac{n^2(x^2+y^2)}{2}+\frac{(1-\mu)q_1}{r_1}+\frac{\mu}{r_2}\nonumber\\&&+\frac{\mu
 A_2}{2r_2^3}+\frac{M_b}{\left(r^2+T^2\right)^{1/2}}\label{eq:OmegaFF}
 \end{eqnarray}
Then the perturbed mean motion $n$ of the primaries is changed into the form  $n^{2}=1+\frac{3A_{2}}{2}+\frac{2M_b}{\left(r_c^2+T^2\right)^{3/2}}$, where $r_c^2=(1-\mu)q_1^{2/3}+\mu^2$, we set $r=r_c=0.08, T=0.01, \mu=0.025$ for further numerical results.
\begin{figure}
\plotone{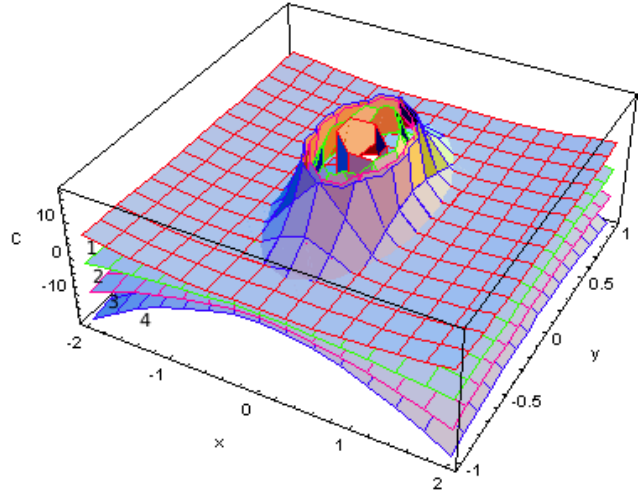}
\caption{The zero velocity curves (1):$M_b=0$, (2):$M_b=0.2$,(3):$M_b=0.4$, (4):$M_b=0.6$, when $\mu=0.025, r=0.8, T=0.01, A_2=0, q=1$} \label{fig:c}
 \end{figure}
The energy integral of the problem is given by  $C=2\Omega-{\dot{x}}^2-{\dot{y}}^2$, where the quantity $C$ is the Jacobi's  constant. The zero velocity curves[see figure  (~\ref{fig:c})]  are given by:
\begin{equation}
C=2\Omega(x,y)\label{eq:C}
\end{equation}
\section{Position of  Equilibrium Points}
\label{sec:existEP}
The  position equilibrium points of Chermnykh's probelm is given by putting  $\Omega_x=\Omega_y=0$ i.e.,  \begin{eqnarray}
&& n^{2}x-\frac{(1-\mu)q_1(x+\mu)}{r^3_1}-\frac{\mu(x+\mu-1)}{r^3_2}\nonumber\\&&-\frac{3}{2}\frac{\mu{A_2}(x+\mu-1)}{r^5_2}-\frac{M_b x}{\left(r^2+T^2\right)^{3/2}}=0\label{eq:eq1pts},\\
&&n^{2}y
-\frac{(1-\mu)q_{1}{y}}{r^3_1}-\frac{\mu{y}}{r^3_2}-\frac{3}{2}\frac{\mu{A_2}y}{r^5_2}\nonumber\\&&-\frac{M_b y}{\left(r^2+T^2\right)^{3/2}} =0\label{eq:eq2pts}\end{eqnarray}
From  equation (~\ref{eq:eq2pts})
\(
n^{2}
-\frac{(1-\mu)q_{1}}{r^3_1}-\frac{\mu}{r^3_2}-\frac{3}{2}\frac{\mu{A_2}}{r^5_2}-\frac{M_b}{\left(r^2+T^2\right)^{3/2}}=0 
\) or $y=0$. 
\begin{figure}
\plotone{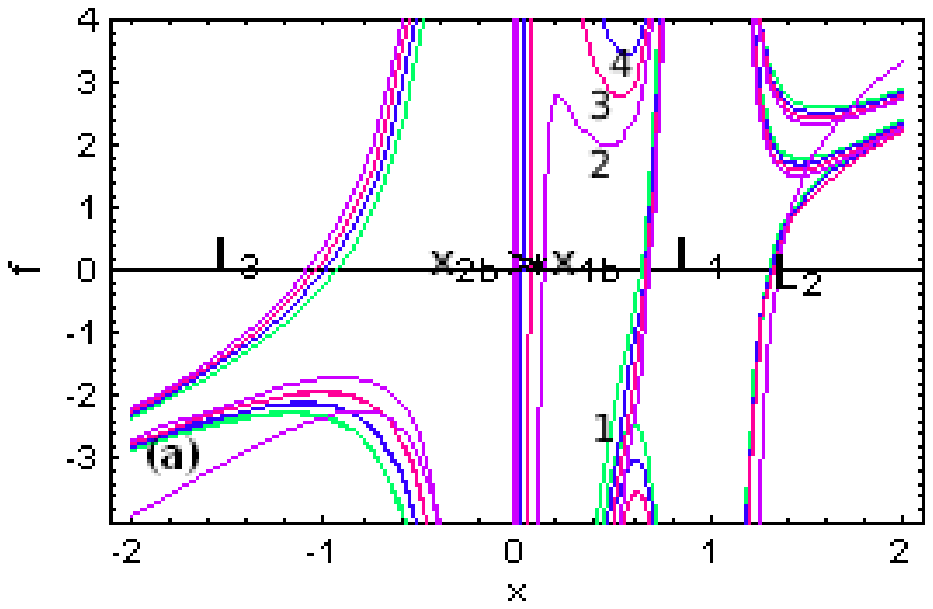}\plotone{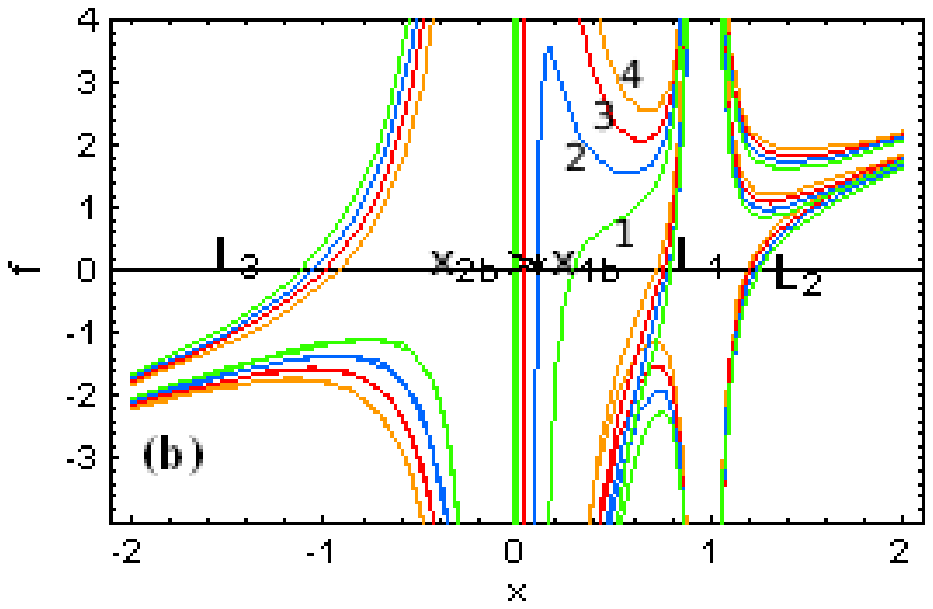}
\caption{Figure show the position of the  collinear  equilibrium points in frame (a) when $q_1=1,A_2=0$ and frame (b) when $q_1=1,A_2=0.02$, the curves are leveled by (1-4) correspond to the mass of belt  $M_b=0, 0.2,0.4,0.6$.}\label{fig:colab}
 \end{figure}
\subsection{Collinear Equilibrium Points}
In this case suppose 
\begin{eqnarray}
&&f(x,y)= n^{2}x-\frac{(1-\mu)q_1(x+\mu)}{r^3_1}-\frac{\mu(x+\mu-1)}{r^3_2}\nonumber\\&&-\frac{3}{2}\frac{\mu{A_2}(x+\mu-1)}{r^5_2}-\frac{M_b x}{\left(r^2+T^2\right)^{3/2}}=0,\label{eq:eq1f},\\
&&g(x,y)=\Biggl[n^{2}
-\frac{(1-\mu)q_{1}}{r^3_1}-\frac{\mu}{r^3_2}-\frac{3}{2}\frac{\mu{A_2}}{r^5_2}\nonumber\\&&-\frac{M_b}{\left(r^2+T^2\right)^{3/2}}\Biggr]y=0\label{eq:eq2g}\end{eqnarray} 
\begin{equation}f(x,0)=P(x)+Q(x)\label{eq:fx0} \end{equation} where 
\begin{eqnarray}
&&P(x)=n^{2}x-\frac{(1-\mu)q_1(x+\mu)}{|x+\mu|^3}-\frac{\mu(x+\mu-1)}{|x+\mu-1|^3}\nonumber\\&&-\frac{3}{2}\frac{\mu{A_2}(x+\mu-1)}{|x+\mu-1|^5}\label{eq:P}
\end{eqnarray}
\begin{equation}
 Q(x)=-\frac{M_b x}{\left(x^2+T^2\right)^{3/2}}=0,\label{eq:Q}
\end{equation}
To investigate the  position of collinear equilibrium points divide the orbital plane $Oxy$ into three parts with respect to the primaries $x\leq-\mu$, $1-\mu\leq x$ and $-\mu<x<1-\mu$, for each part the function $P(x)$ is defind as follows: 
\begin{equation}
P(x)=\begin{cases} n^{2}x+\frac{(1-\mu)q_1}{(x+\mu)^2}+\frac{\mu}{(x+\mu-1)^2}\\\ +\frac{3}{2}\frac{\mu{A_2}}{(x+\mu-1)^4} & \text{If}\  x<-\mu,\\\
n^{2}x-\frac{(1-\mu)q_1}{(x+\mu)^2}-\frac{\mu}{(x+\mu-1)^2}\\-\frac{3}{2}\frac{\mu{A_2}}{(x+\mu-1)^4} & \text{If}\ 1-\mu<x,\\
n^{2}x-\frac{(1-\mu)q_1}{(x+\mu)^2}+\frac{\mu}{|x+\mu-1|^3}\\+\frac{3}{2}\frac{\mu{A_2}}{(x+\mu-1)^4} & \text{If}\  -\mu <x<1-\mu
 \end{cases}
\end{equation}

When $x\in(-\infty,-\mu)$,  $\lim_{x\rightarrow -\infty}f(x,y) < 0$, $\lim_{x\rightarrow -\mu^-}f(x,y)>0$, $P'(x)>0, Q'(x)<0$ so there is a point $x_3$ for which $f(x_3,0)=0$. If  $x\in(1-\mu, \infty)$, $\lim_{x\rightarrow (1-\mu)^+}f(x,y)<0$,   $\lim_{x\rightarrow \infty}f(x,y)>0$,\ $P'(x)>0, Q'(x)>0$ so there is a point $x_2$ for which $f(x_2,0)=0$. Now $x\in(-\mu, 1-\mu)$, $\lim_{x\rightarrow 0^+}f(x,y)>0$,\ $\lim_{x\rightarrow(1- \mu)^-}f(x,y)>0$ this implies that an even(or zero) number real roots of $f(x,y)$ exists in this range. If  $x\in (-\mu,0)$, consider two cases (i) $x\in (-\mu,-\frac{T}{\sqrt{2}})$ and (ii) $x\in (-\frac{T}{\sqrt{2}}), 0)$,  we obtained $Q(-\frac{T}{\sqrt{2}})+P(-\frac{T}{\sqrt{2}})>0$ i.e. f$(x,0)\neq0$ so there is no equilibrium point at $x=-\frac{T}{\sqrt{2}}$. If $Q(-\frac{T}{\sqrt{2}})+P(-\frac{T}{\sqrt{2}})>0$, we obtained    $\lim_{x\rightarrow 0^-}f(x,y)<0$, $\lim_{x\rightarrow- \mu^+}f(x,y)<0$ so there exists two new  equilibrium points $x_{b1}\in(-\frac{T}{\sqrt{2}},0)$ and $x_{b2}\in(-\mu, -\frac{T}{\sqrt{2}})$   for which $f(x,0)=0$. But $Q(-\frac{T}{\sqrt{2}})+P(-\frac{T}{\sqrt{2}})<0$ the  there exists no equilibrium points in $(-\mu,0)$. If $T<\sqrt{2}\mu$ and $Q(-\frac{T}{\sqrt{2}})+P(-\frac{T}{\sqrt{2}})>0$ then we  have  two new equilibrium points. Hence we have found there are five equilibrium points on the $x$-axis for the given system. The position of above points are presented graphically by frames a, b in figure (~\ref{fig:colab}) the curves are leveled by (1-4) correspond to the mass of belt  $M_b=0, 0.2,0.4,0.6$.
\subsection{Triangular Equilibrium Points}
\label{subsec:TrigPts}
The triangular equilibrium points are given by putting  $\Omega_x=\Omega_y=0$, $y\neq{0}$. Using the method as \citet*{kushvah-2008} then from equations (~\ref{eq:Omegax}) and (~\ref{eq:Omegay}) we obtained:
\begin{eqnarray}
r_1&&=q_1^{1/3}\left[1-\frac{A_2}{2}\right.\nonumber\\&&\left.+\frac{(1-2r_c)M_b \left(1-\frac{3\mu A_2}{2(1-\mu)}\right)}{3\left(r_c^2+T^2\right)^{3/2}}\right]\label{eq:r1},\\ 
r_2&&=1+\frac{\mu(1-2r_c)M_b}{3\left(r_c^2+T^2\right)^{3/2}}\label{eq:r2}
\end{eqnarray}
From above, the triangular equilibriumpoints are as:
\begin{eqnarray}
&&x=-\mu+\frac{q_1^{2/3}}{2}\biggl[1-\frac{A_2}\nonumber\\&&+\frac{\mu(2r_c-1)M_b A_2}{(1-\mu)\left(r_c^2+T^2\right)^{3/2}}\biggr]\label{eq:xl4}\\
 &&y=\pm\frac{q_1^{2/3}}{2}\biggl[\left(4-q_1^{2/3}\right)+2\left(q_1^{2/3}-2\right)A_2\nonumber\\&&+\frac{4(2r_c-1)M_b\left[\left\{\left(q_1^{2/3}-3\right)-\frac{3\mu A_2\left(q_1^{2/3}-3\right)}{2(1-\mu)}\right\}\right]}{3\left(r_c^2+T^2\right)^{3/2}}\biggr]^{1/2}\label{eq:yl4} 
 \end{eqnarray}
All these results are similar with \citet*{Szebehely1967}, \citet*{Ragosetal1995},  \citet*{JiangYeh2006Ap&SS} and others. 
\section{Linear Stability}
\label{sec:lstb_with_PR}
To study the linear stability of any equilibrium point  change the origin of the coordinate system to its position $(x*,y*)$ by means of  $x=x*+\alpha$,\,  $y=y*+\beta$, where  $\alpha=\xi e^{\lambda{t}}$,\ $\beta=\eta e^{\lambda{t}}$ are the small displacements  $\xi,\eta$,\  $\lambda$ these parameters, have to be determined.  Therefore the equations of perturbed motion corresponding to the system of equations (~\ref{eq:Omegax}), (~\ref{eq:Omegay}) may be written as follows:
\begin{align}
\ddot{\alpha}-2n\dot{\beta} &= {\alpha}{\Omega^*_{xx}}+{\beta}{\Omega^*_{xy}} \\
\ddot{\beta}+2n\dot{\alpha}&= {\alpha}{\Omega^*_{yx}}+{\beta}{\Omega^*_{yy}}
\end{align}
where superfix $*$ is corresponding to the equilibrium points.
 \begin{align}(\lambda^2-{\Omega^*_{xx}})\xi
+(-2n\lambda-{\Omega^*_{xy}})\eta&=0\label{eq:lambda_x}\\
(2n\lambda-{\Omega^*_{yx}})\xi
+(\lambda^2-{\Omega^*_{yy}})\eta&=0\label{eq:lambda_y}
\end{align}
Now  above system has singular solution if,
\[
\begin{vmatrix}
\lambda^2-\Omega^*_{xx}& -2n\lambda-\Omega^*_{xy} \\2n\lambda-\Omega^*_{yx}& \lambda^2-\Omega^*_{yy}\\
\end{vmatrix}
=0
\]
\begin{equation}
\Rightarrow \quad \lambda^4+b\lambda^2+d=0 \label{eq:cheq}
\end{equation}
At the equilibrium points equations (~\ref{eq:ux}),(~\ref{eq:uy}) gives us the following:
\begin{eqnarray*}
 &&b=2n^2-f_*-\frac{3\mu{A_2}}{{r^5_2}_*}+\frac{3M_b T^2}{\left(r^2_*+T^2\right)^{5/2}}\\
&&d=9\mu\amc{y^2_*}\left[\frac{q_1}{r^5_{1*}r^5_{2*}}\right.\\&&\left.+\frac{3M_b}{\left(r^2_*+T^2\right)^{5/2}}\left\{\frac{\mu q_1}{{r_1^5}_*}+\frac{(1-\mu)\zabs}{r_{2*}^5}\right\}\right]
 \end{eqnarray*}
where \(f_*=\frac{(1-\mu)q_1}{r^3_{1_*}}+\frac{\mu}{r^3_{2*}}\left(1+\frac{3}{2}\frac{A_2}{r^2_{2*}}\right)+\frac{3M_b }{\left(r^2_*+T^2\right)^{5/2}}\).

\begin{figure}[h]
\plotone{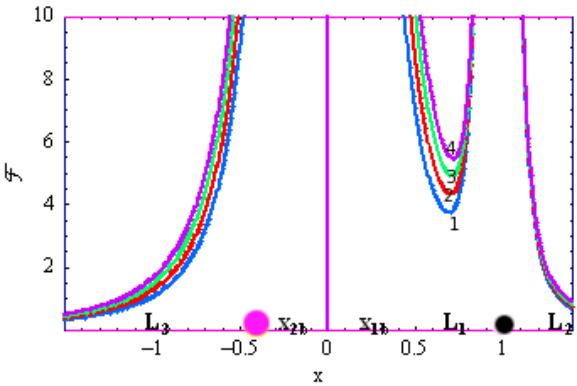} \caption{$x-f$, curves correspond to (1):$M_b=0$, (2):$M_b=0.2$,(3):$M_b=0.4$, (4):$M_b=0.6$, when $\mu=0.025,r=0.8, T=0.01$, $ A_2=0$ and  $q_1=1$}\label{fig:col_f}
 \end{figure}

 The points $L_1,L_2,L_3$, $x_{b1},x_{b2}$   lie along the line joining the primaries, from (~\ref{eq:r1}) we have  $r_1\approx \frac{q_1^{1/3}}{n^2}$. In this case 
\begin{eqnarray}f&=&\frac{(1-\mu)q_1}{|x+\mu|^3}+\frac{\mu}{|x+\mu-1|^3}\left(1+\frac{3}{2}\frac{A_2}{|x+\mu-1|^2}\right)\nonumber\\&&+\frac{3M_b }{\left(x^2+T^2\right)^{5/2}}\label{eq:fcol}
\end{eqnarray}
With the help of equation (~\ref{eq:fcol}) we  have drawn figure  ~\ref{fig:col_f}, the different  curves (1)-(4) correspond to $M_b=0.0, 0.2, 0.4, 0.6$. From this figure we observe that  for each  $L_1,L_2,L_3$, $x_{b1}$ and $x_{b2}$,  $f>1$, in this case  characteristic equation (~\ref{eq:cheq}) has at least one positive root, this implies that collinear equilibrium  points are unstable in the sense of Lyapunov. 

\begin{figure}
\plotone{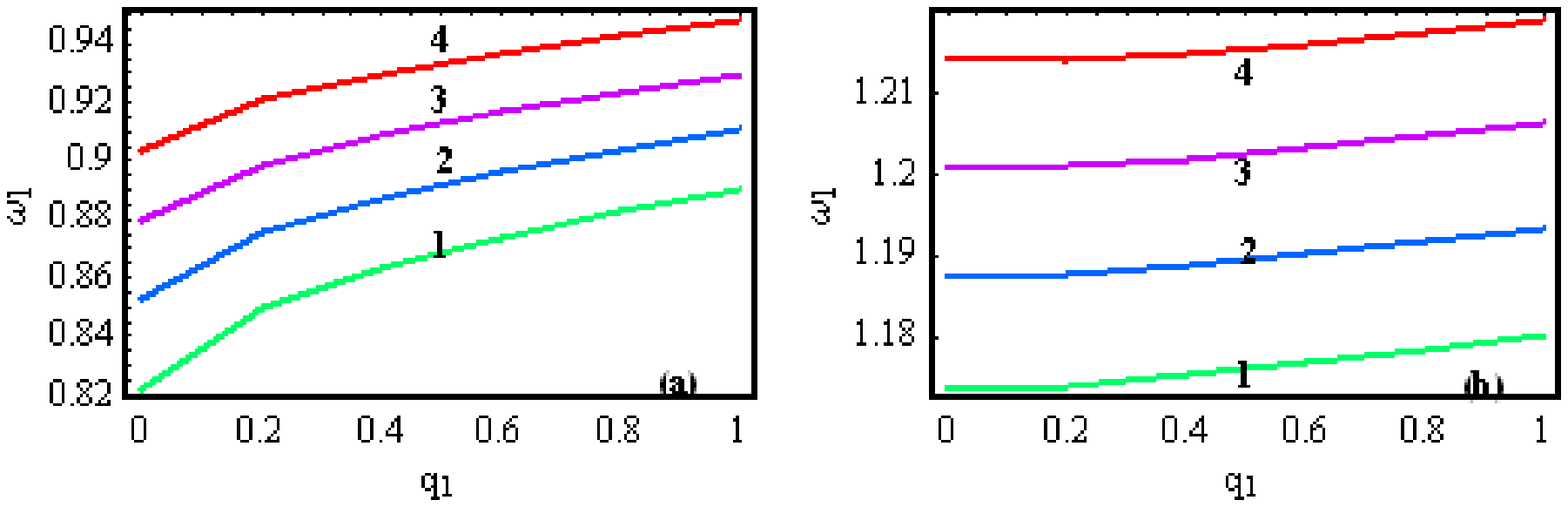}\plotone{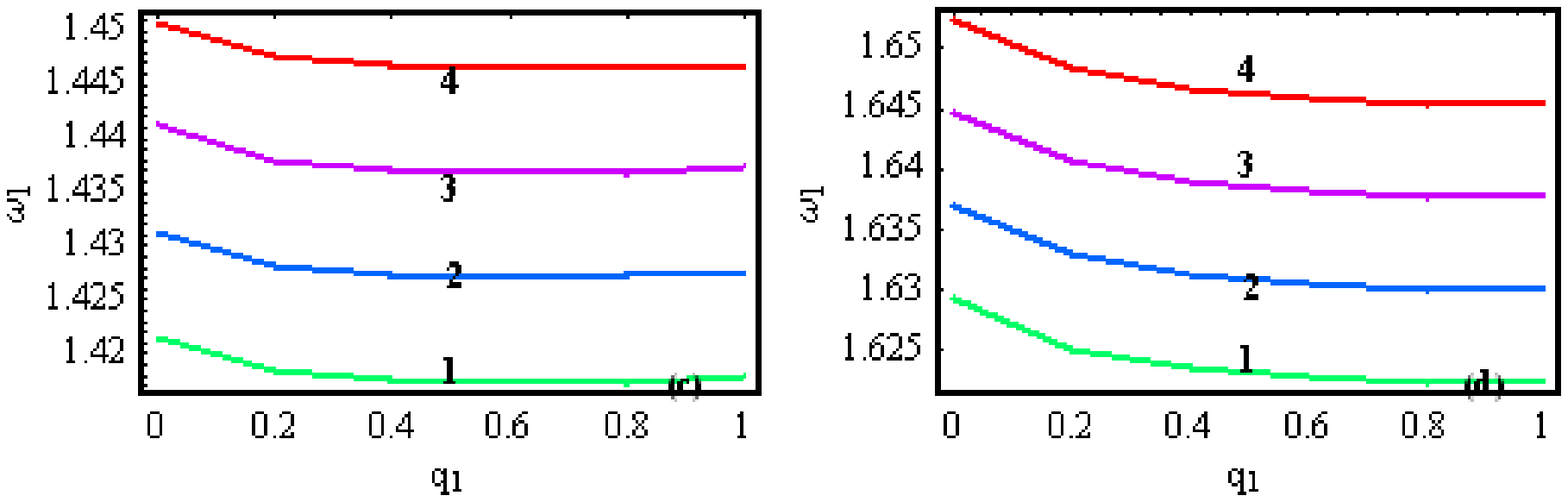}
\caption{$\omega_1-q_1$,frames correspond to (a):$M_b=0$, (b):$M_b=0.2$,(c):$M_b=0.4$, (d):$M_b=0.6$, curves $(1):A_2=0$, $(2):A_2=0.02$, $(3):A_2=0.04$, $(4):A_2=0.06$, $\mu=0.025,r=0.8, T=0.01 $}\label{fig:omega1A}
 \end{figure}

\begin{figure}
\plotone{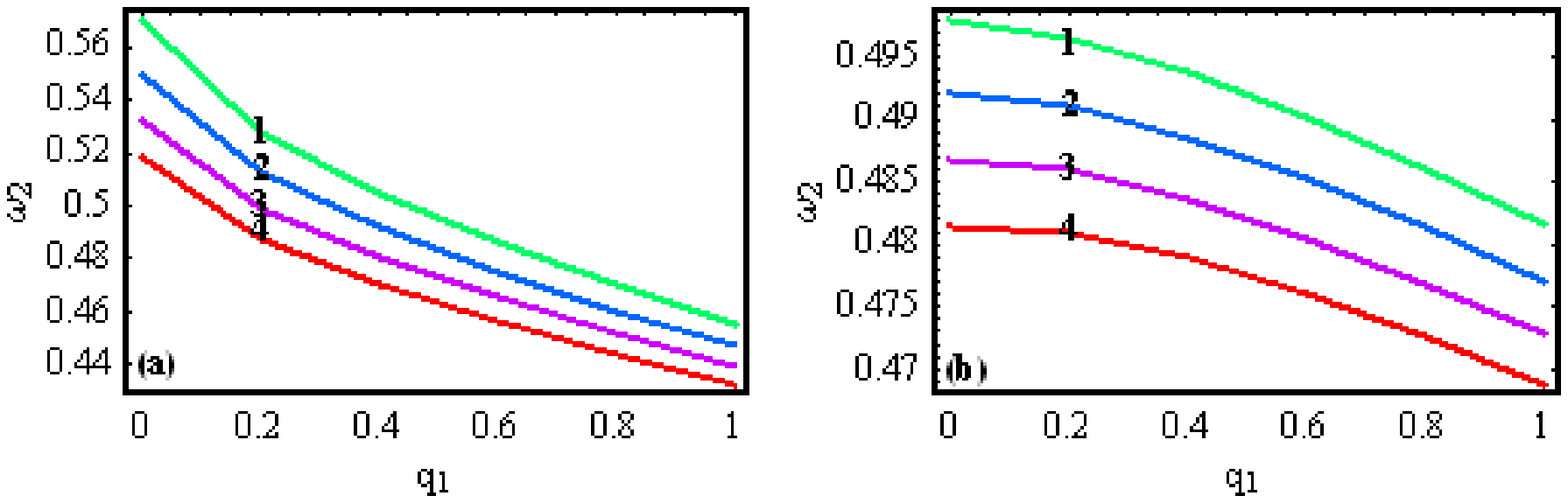}\plotone{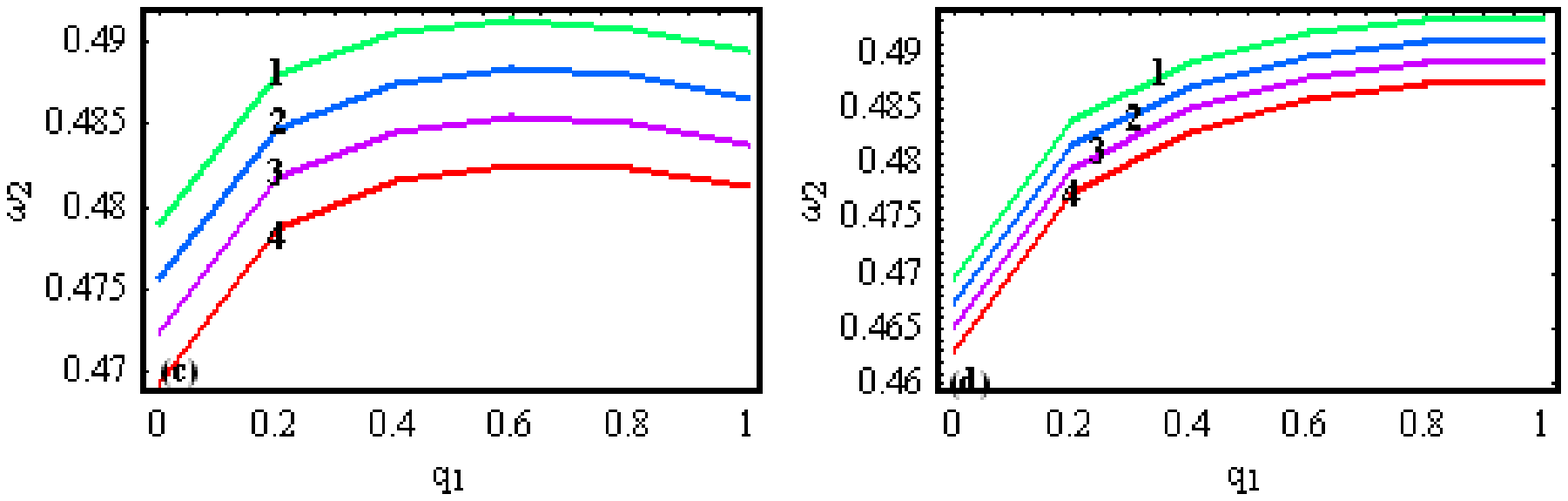}
\caption{$\omega_2-q_1$, frames correspond to (a):$M_b=0$, (b):$M_b=0.2$,(c):$M_b=0.4$, (d):$M_b=0.6$, curves $(1):A_2=0$, $(2):A_2=0.02$, $(3):A_2=0.04$, $(4):A_2=0.06$, $\mu=0.025,r=0.8, T=0.01 $}\label{fig:omega2A}
 \end{figure}

\begin{figure}
\plotone{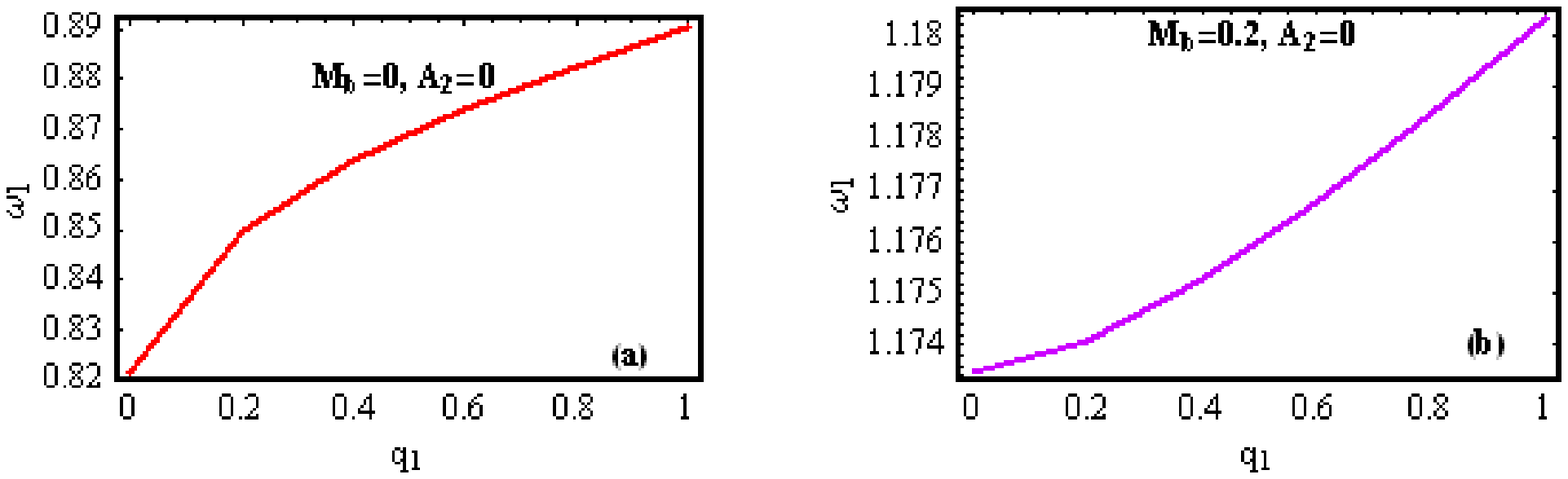}\plotone{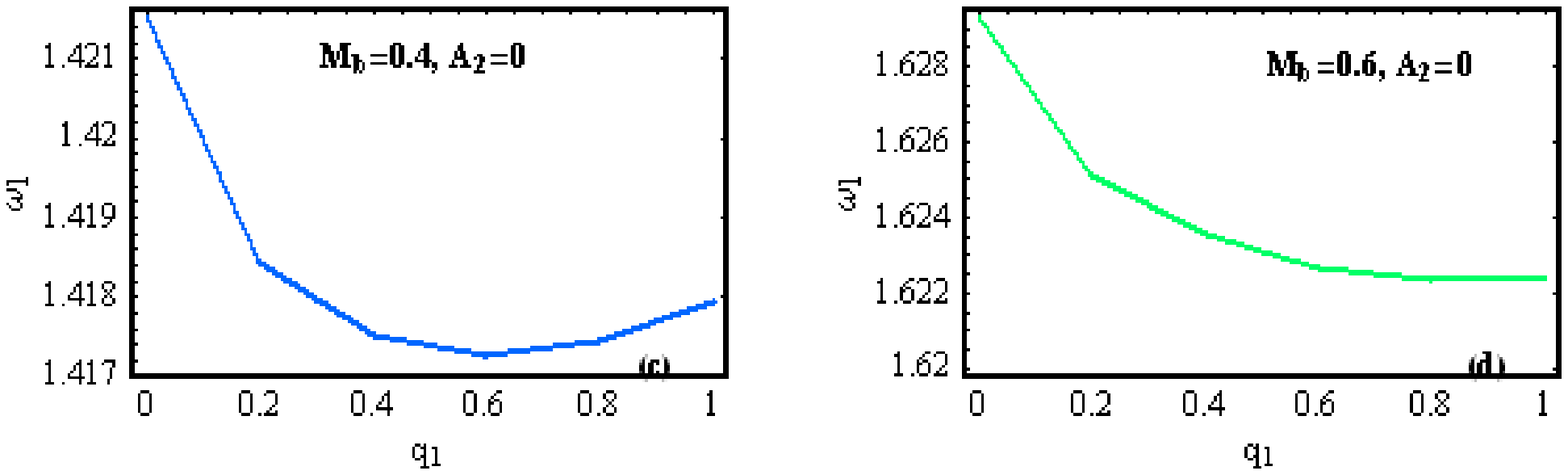}
\caption{$\omega_1-q_1$ frames (a):$M_b=0$, (b):$M_b=0.2$,(c):$M_b=0.4$, (d):$M_b=0.6$, when $\mu=0.025, r=0.8, T=0.01$}\label{fig:omega1}
 \end{figure}

\begin{figure}
\plotone{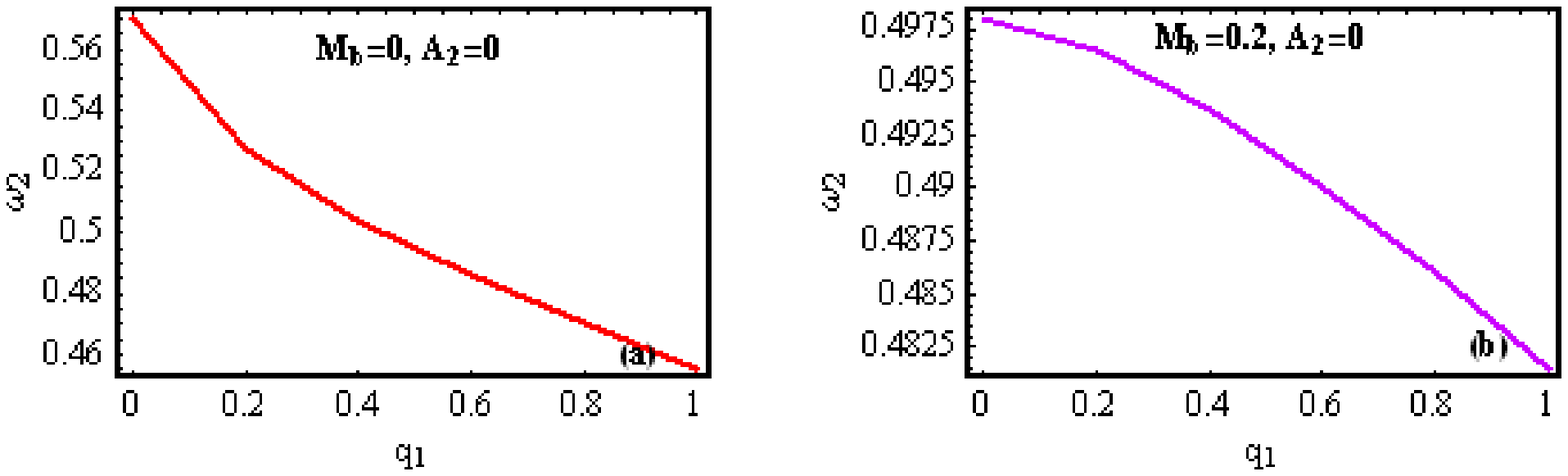}\plotone{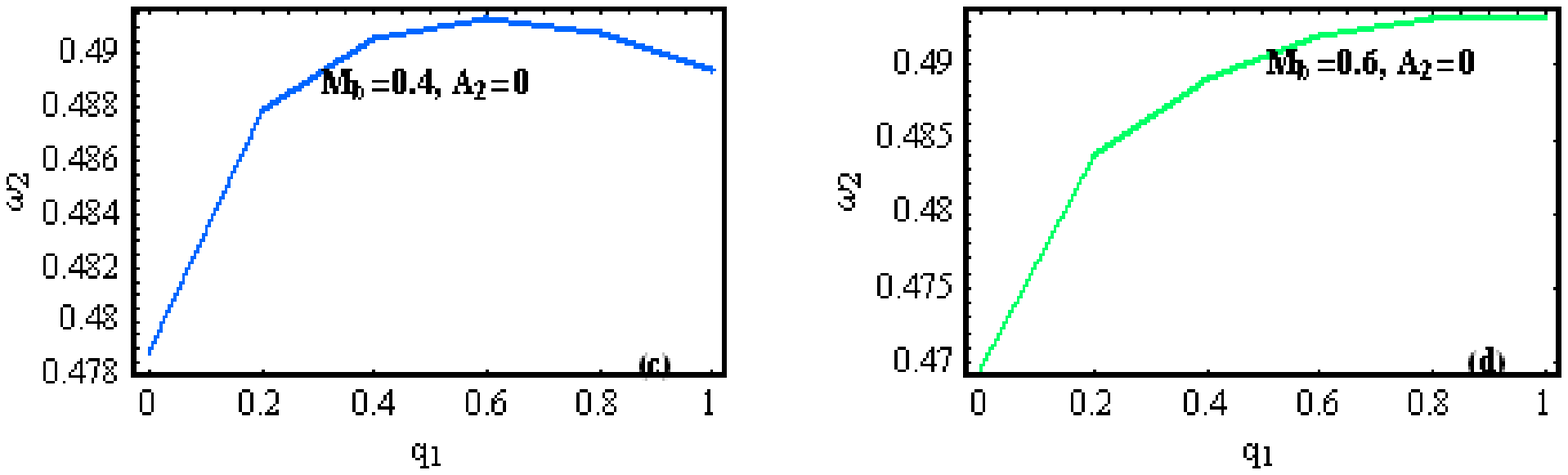}
\caption{$\omega_2-q_1$, frames (a):$M_b=0$, (b):$M_b=0.2$,(c):$M_b=0.4$, (d):$M_b=0.6$, when $\mu=0.025, r=0.8, T=0.01$}\label{fig:omega2}
 \end{figure}

Now we have to study the linear stability of triangular equilibrium points, in this regard  we obtained \(f=n^2,\  b=9\mu\amc g\), where \(g={y^2_*}\left[\frac{q_1}{r^5_{1*}r^5_{2*}}\right.\\\left.+\frac{3M_b}{\left(r^2_*+T^2\right)^{5/2}}\left\{\frac{\mu q_1}{{r_1^5}_*}+\frac{(1-\mu)\zabs}{r_{2*}^5}\right\}\right]\).  From characteristic equation (~\ref{eq:cheq}) we obtained:
\begin{equation}
 \lambda^2=\frac{-b\pm{(b^2-4d)}^{1/2}}{2}\label{eq:cheqNpr}
\end{equation}
For stable motion  $0<4d<b^2$, i.e.
\begin{equation*}
{\bigl(n^2-3\mu{A_2}\bigr)^2}>36\mu\amc g
\end{equation*}
In classical case $A_2=0,M_b=0$, $q_1=1$, $n=1$, we have following:
\(
1>27\mu\amc\  \Rightarrow \quad  \mu <0.0385201.
\)
Using equation (~\ref{eq:cheqNpr}) we obtained  imaginary roots $\lambda_{1,2}=\pm \mathbf{i}\omega_1$, $\lambda_{3,4}=\pm \mathbf{i}\omega_2$, $\mathbf{i}=\sqrt{-1}$.
The characteristic frequencies $\omega_{1,2}(0<\omega_2<\omega_1)$  are presented by frames (a)-(d) of  (~\ref{fig:omega1A}--~\ref{fig:omega2}) in parameter plots $\omega_i-q_1$ for different values of $A_2, q_1, M_b$. They are given in  table(~\ref{tbl-2}). We observe that they are decreasing function of radiation pressure and increasing functions of  $A_2, M_b$.   Hence the triangular  equilibrium points are stable in the sense of Lyapunov stability provided  $\mu< \mu_{Routh}=0.0385201$. 
In this case we obtained the three main cases of resonances:
\begin{equation}
\omega_1-k\omega_2=0,\quad k=1,2,3\label{eq:Reson}
\end{equation}
For $k=1$ we have positive stable resonance and for $k=2,3$ we have unstable resonances. Using (~\ref{eq:cheqNpr})  and (~\ref{eq:Reson}) we obtained a root of mass parameter:
\begin{eqnarray}
\mu_k=\frac{3g+2 K b_1b_2-\sqrt{g}\sqrt{9g-4Kb_1^2+12b_1b_2}}{6(g+K b_2^2)}
\end{eqnarray}
where $K=\frac{k^2}{(k^2+1)^2}$,\ \(b_1=n^2+\frac{2rM_b }{\left(r^2+T^2\right)^{3/2}}+\frac{3M_b T^2}{\left(r^2+T^2\right)^{5/2}}\), \(b_2=A_2\left[1+\frac{5(2r-1)M_b }{\left(r^2+T^2\right)^{3/2}} \right]\).
Now we suppose  $q_1=1-\epsilon$, with $|\epsilon|<<1$, neglecting higher  order terms, we obtained the critical mass parameter values corresponding to $k=1,2,3$ as :
\begin{eqnarray}
\mu_1&=&0.0385208965 + 0.0375419787A_2\nonumber\\&&- 
      0.0089174706\epsilon- 
      0.0678734040 M_b\\
\mu_2&=& 0.0242938971+ 0.0254350205 A_2\nonumber\\&&- 
    0.0055364958\epsilon- 
      0.0421398438 M_b\\
\mu_3&=&0.0135160160+ 0.0148764140A_2\nonumber\\&&- 
    0.0030452832\epsilon- 
      0.0231785159 M_b
\end{eqnarray}
\begin{figure}
\plotone{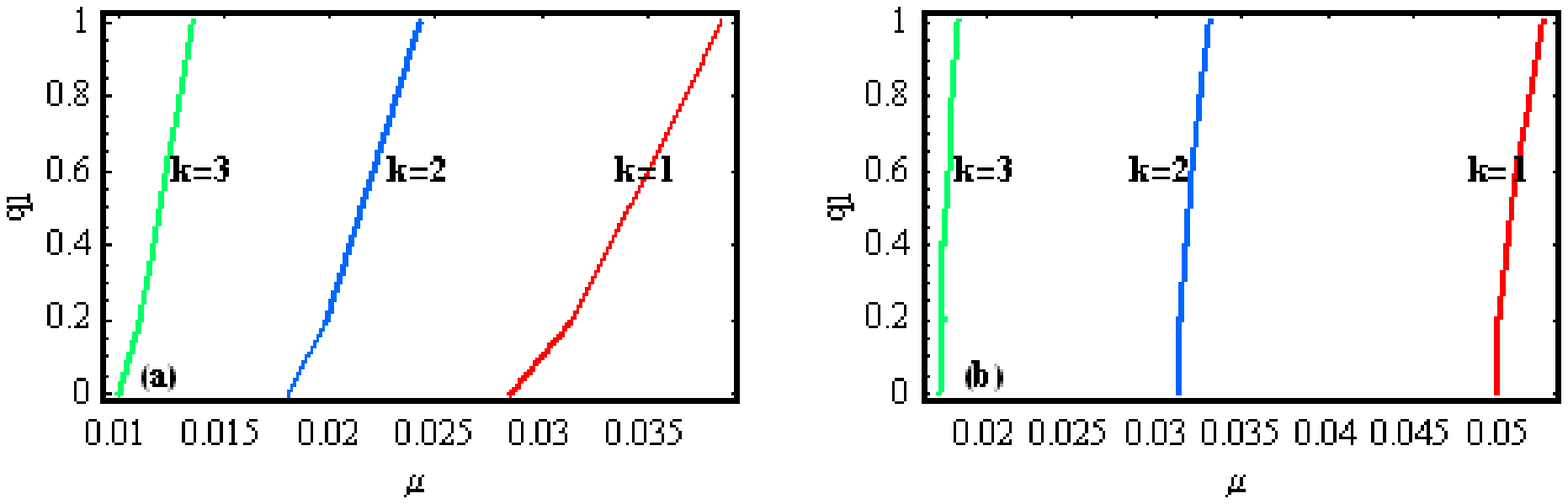}\plotone{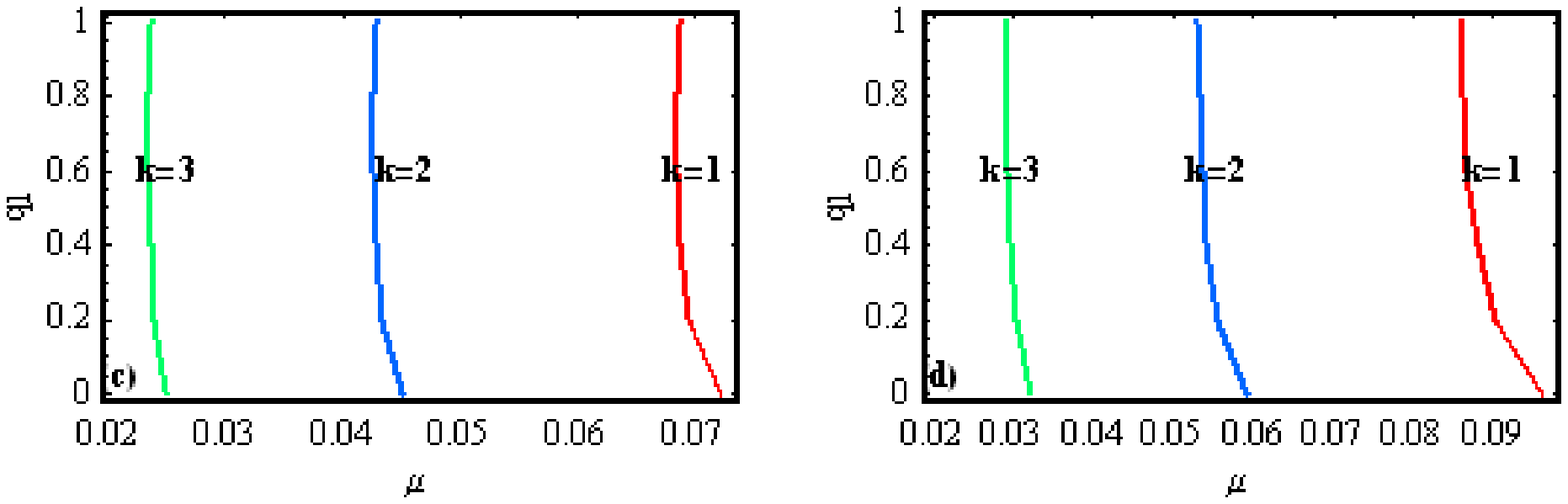}
\caption{$\mu-q_1$  frames correspond to (a):$M_b=0$, (b):$M_b=0.2$,(c):$M_b=0.4$, (d):$M_b=0.6$, when $A_2=0, k=1-3$, $r=0.8, T=0.01$}\label{fig:mukA0}
 \end{figure} 
\begin{figure}
\plotone{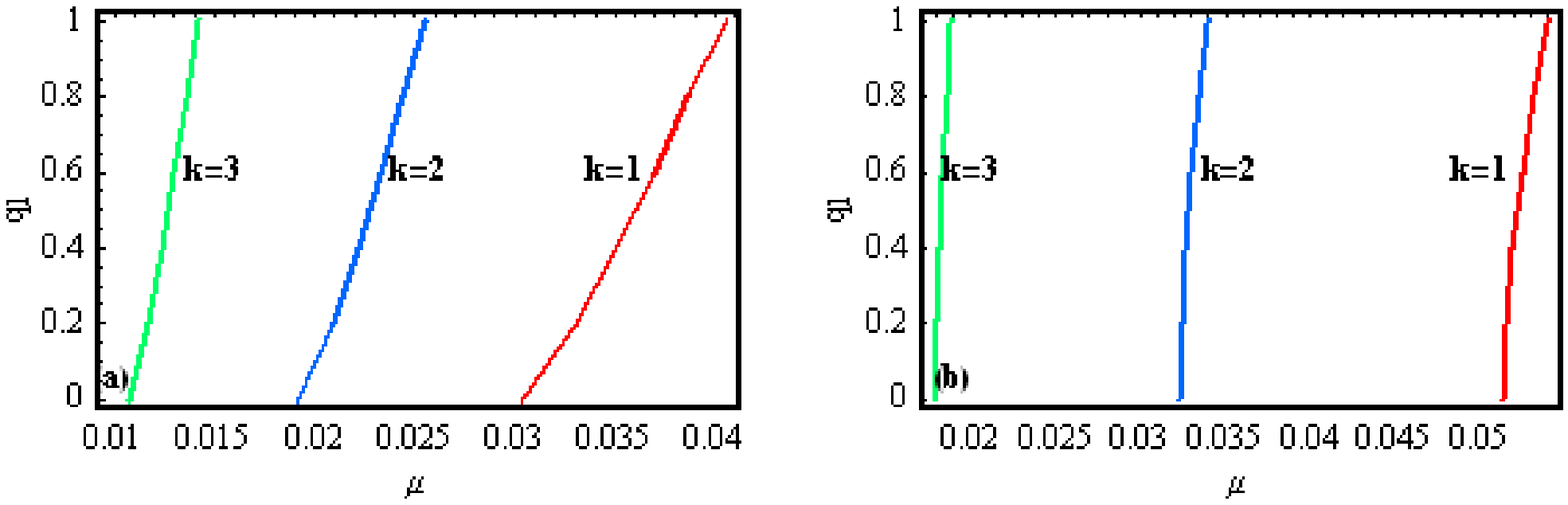}\plotone{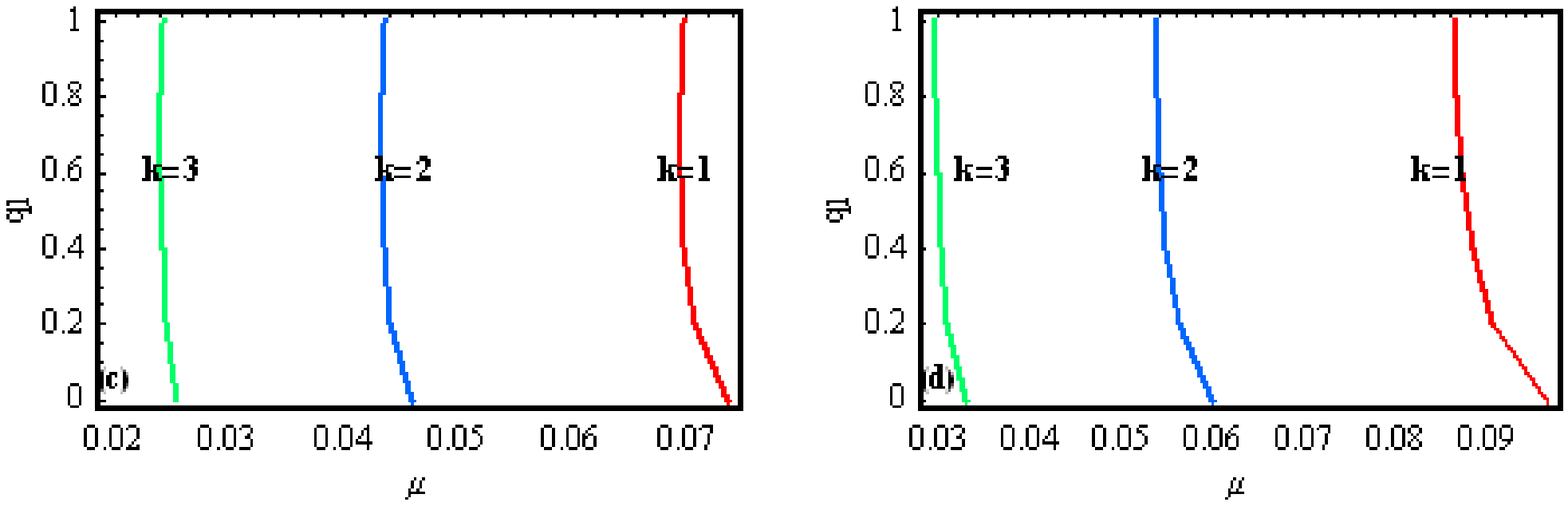}
\caption{$\mu-q_1$ frames correspond to (a):$M_b=0$, (b):$M_b=0.2$,(c):$M_b=0.4$, (d):$M_b=0.6$, when $A_2=0.02, k=1-3$, $r=0.8, T=0.01$}\label{fig:mukA2}
 \end{figure} 
 The linear stability region and   main resonance curves $k=1,2,3$ are shown by $\mu-q_1$ parameter space frames (a)-(d) in figures (~\ref{fig:mukA0},~\ref{fig:mukA2}). The  curve corresponding to $k=1,(q_1=1, A_2=0,M_b=0, \mu_1=\mu_{Routh}=0.038521)$ is actual boundary of the stability region. The critical values of mass parameter $\mu_k$ are presented in  table (~\ref{tbl-3}) for various values of $q_1, A_2, M_b$. The classical critical values of $\mu$ are similar to \citet*{Deprit1967}. These results are similar to the results of  \citet*{Markellos1996Ap&SS}, \citet*{kushvah-2008} and others.
 We observe that the effect of radiation pressure reduces the  linear stability zones and the  $\mu_k$ is an increasing function of $M_b$, $A_2$.
\section{Conclusion}
\label{sec:con}
 The points $L_1,L_2,L_3$, $x_{b1},x_{b2}$  lie along the line joining the primaries. We observe that the effect of radiation pressure reduces the  linear stability zones,  these are also affected by belt and the oblateness of second primary. The collinear equilibrium points  are unstable while triangular  equilibrium points are stable in the sense of Lyapunov stability provided  $\mu< \mu_{Routh}=0.0385201$.

\acknowledgments{I am very thankful  to Dr. Uday Dolas, Dr. Deepak Singh  for their persuasion. I am also thankful to Mrs. Snehlata Kushwah for her loving support. I am grateful to IUCAA Pune for financial assistance to visit library and computer facilities.}
\bibliographystyle{spbasic} 

\clearpage
\begin{deluxetable}{rrrrrrrrrr}
\tabletypesize{\scriptsize}
\rotate
\tablecaption{$\omega_{1,2}$ when  $r=0.8,T=0.01, \mu=0.025$\label{tbl-2}}
\tablewidth{0pt}
\tablehead{
\colhead{$A_2$}& \colhead{$q_1$} & \colhead{$\omega_1:M_b=0$} & 			\colhead{$\omega_2:M_b=0$} & \colhead{$\omega_1:M_b=0.2$} & \colhead{$\omega_2:M_b=0.2$} &
\colhead{$\omega_1:M_b=0.4$} & \colhead{$\omega_2:M_b=0.4$}&
\colhead{$\omega_1:M_b=0.6$} & \colhead{$\omega_2:M_b=0.6$}
}
\startdata
0.0 & 1.0 & 0.890141&0.455686&1.18033&0.4815237&1.41795&0.489382&1.62232&0.493127\\
&0.75&0.880622&0.47382&1.17804&0.487083&1.41737&0.491041&1.62238&0.492928\\
&0.5&0.869076&0.494679&1.17602&0.491953&1.41733&0.491159&1.62304&0.490779\\
&0.25&0.853749&0.520684&1.17436&0.495895&1.41812&0.488893&1.62461&0.485532\\
&0.0& 0.821584&0.570088&1.17349&0.497955&1.4215&0.478958&1.62926&0.4697\\
&&&&&&&&&\\
0.02& 1.0& 0.910283&0.447086&1.1934&0.477057&1.42768&0.486598&1.63005&0.491211\\
&0.75&0.901845&0.463871&1.19127&0.482336&1.42716&0.48813&1.63013&0.490934\\
&0.5&0.891743&0.483005&1.18941&0.486917&1.42716&0.488124&1.6308&0.488709\\
&0.25&0.878605&0.506511&1.18791&0.490552&1.42798&0.485737&1.63238&0.48339\\
&0.0&0.852388&0.549485&1.18724&	1.43137&1.43137&0.475658&1.63701&0.467476\\
&&&&&&&&&\\
0.04& 1.0 & 0.929538&0.439272&1.20624&0.472778&1.43732&0.483883&1.63772&0.489328\\
&0.75&0.921985&0.45491&1.20426&0.477788&1.43685&0.48529&1.63783&0.488972\\
&0.5&0.913034&0.47262&1.20254&0.482095&1.43689&0.485164&1.63851&0.486674\\
&0.25&0.901559&0.494157&1.2012&0.485439&1.43774&0.482658&1.6401&0.481283\\
&0.0&0.879387&0.532615&1.2007&0.486672&1.44113&0.472436&1.64471&0.465288\\
\enddata
\tablecomments{Table \ref{tbl-2}  presents the roots of characteristic equation(~\ref{eq:cheq}).}
\end{deluxetable}
\begin{deluxetable}{rrrrrrrrrr}
\tabletypesize{\scriptsize}
\rotate
\tablecaption{$\mu_k(A_2,M_b)$, when $ r=0.8,T=0.01$\label{tbl-3}}
\tablewidth{0pt}
\tablehead{
\colhead{$q_1$}&\colhead{$k$} & \colhead{$\mu_k(0,0)$} & 			\colhead{$\mu_k(0,0.02)$} & \colhead{$\mu_k(0,0.4)$} & \colhead{$\mu_k(0,0.6)$} &
\colhead{$\mu_k(0.02,0)$} & \colhead{$\mu_k(0.02,0.2)$}&\colhead{$\mu_k(0.02,0.4)$}&\colhead{$\mu_k(0.02,0.6)$}
}
\startdata
1.0&1&0.0385209&0.0525812&0.0688051&0.0861218&0.0404877&0.0539744&0.0696964&0.0863953\\
&2&0.0242939&0.0329695&0.0428408&0.0532015&0.0255597&0.0339343&0.0435889&0.0537117\\
&3&0.013516&0.0182676&0.0236236&0.0291855& 0.0142309&0.0188392&0.0241121&0.0295953\\
&4&0.00827037&0.0111565&0.014396&0.0177443&0.00871096&0.0115163&0.0147154&0.01803\\
&5&0.0055092&0.00742441&0.00956953&0.0117815&0.0058038&0.00766762&0.00978937&0.0119837\\
&&&&&&&&&\\
0.75&1&0.0363201&0.051579&0.0684542&0.0863238&0.0382382&0.0529886&0.0693753&0.0866222\\
&2&0.0229262&0.0323547&0.0426289&0.0533212&0.024159&0.0333263&0.0433931&0.0538481\\
&3&0.0127632&0.0179323&0.0235093&0.0292495&0.0134588&0.0113141&0.0240058&0.0296688\\
&4&0.0078121&0.0109532&0.014327&0.0177827&0.00824062&0.0113141&0.014651&0.0180743\\
&5&0.00520474&0.00728962&0.00952388&0.0118069&0.00549121&0.0075334&0.00974673&0.012013\\
&&&&&&&&&\\
0.5&1&0.0341355&0.0507482&0.0685365&0.0872546&0.0359977&0.0521785&0.0694903&0.0875708\\
&2&0.0215661&0.0318447&0.0426786&0.0538727&0.0227616&0.0328263&0.0434633&0.054418\\
&3&0.0120136&0.0176539&0.0235361&0.0295437&0.0126877&0.0182322&0.0240439&0.0299757\\
&4&0.00735548&0.0107843&0.0143431&0.0179594&0.00777057&0.0111476&0.0146741&0.0182594\\
&5&0.00490128&0.00717768&0.00953459&0.0119234&0.00517872&0.00742292&0.00976201&0.0121354\\
\enddata
\tablecomments{Table \ref{tbl-3} is presents the roots of characteristic equation(~\ref{eq:cheq})}
\end{deluxetable}
\clearpage
\end{document}